\documentclass[journal]{IEEEtran}
\usepackage{cite}
\usepackage{hyperref}

\ifCLASSINFOpdf
\else
\fi
\usepackage{graphicx}
\usepackage{epstopdf}

\usepackage{amsmath}
\usepackage{amssymb}

\usepackage{multirow}

\usepackage{array}
\usepackage{float}
\usepackage{cancel}
\usepackage{graphicx}

\usepackage{url}

\usepackage[hyphenbreaks]{breakurl}

\newcommand{\matpower}{M{\sc atpower}}

\begin{document}
%

%
%
\title{AC Power Flow Data in MATPOWER and QCQP format: iTesla, RTE Snapshots, and PEGASE}
%
%

%
%
%

\author{C\'edric~Josz,
        St\'ephane~Fliscounakis,
        Jean~Maeght,
        and~Patrick~Panciatici
\thanks{The authors are with the French transmission system operator RTE, 9 rue de la Porte de Buc, BP 561,
F-78000 Versailles, France. E-mails:
\texttt{firstname.lastname@rte-france.com}
}}

\maketitle

\begin{abstract}
In this paper, we publish nine new test cases in \matpower{} format. Four test cases are French very high-voltage grid generated by the offline plateform of iTesla: part of the data was sampled. Four test cases are RTE snapshots of the full French very high-voltage and high-voltage grid that come from French SCADAs via the Convergence software. The ninth and largest test case is a pan-European ficticious data set that stems from the PEGASE project. It complements the four PEGASE test cases that we previously published in \matpower{} version 5.1 in March 2015. We also provide a MATLAB code to transform the data into standard mathematical optimization format. Computational results confirming the validity of the data are presented in this paper.
\end{abstract}


%
\IEEEpeerreviewmaketitle

%
%
\section{Introduction}
%
%

\IEEEPARstart{T}{he} purpose of our grid data publications is to contribute to the progress of the power systems scientific community. As power systems practitioners, we definitely need improved power flow computation methods. By sharing data sets that we use on a daily basis, we hope to help the community develop faster and more reliable optimal power flow methods.

In arXiv, in the download table (other format), source code for the article is delivered as a gzipped tar (.tar.gz) file. It contains the nine test cases as MATLAB .m files.

A lot of European public grid data are already available in various formats~\cite{CIMentsoe,eurodata}. The advantage of the data we are providing is that we have written them in the format of \matpower{}~\cite{matpower}. It also worth to mention the major work made in \cite{nesta2015} to document all of the
AC power transmission system data that is publicly available.

%
%
\section{Origin of the data}
\label{sec:Origin of the data}
%
%

%
\subsection{iTesla}
\label{subsec:iTesla}
%

iTesla stands for: Innovative Tools for Electrical System Security within Large Areas; it is a large collaborative R\&D project co-funded by the EC 7th Framework Programme. Detailed information may be found on the web site of the project~\cite{iTesla} and was presented during dissemination events \cite{dissiTesla,finaliTesla}. In the offline platform of iTesla, a monte-carlo process is run. Loads and uncontrollable generation (mainly solar and wind power) are sampled. From these sampled values, and for each monte-carlo run, a full AC network state is built to serve as a starting point for time domain simulations (e.g. with Eurostag simulation software). Thousands of such processes were run during iTesla project, using High Performance Computing facilities (10,000 cores). Only 4 are published here. 2 of them contain French VHV grid, and 2 of them contain French VHV grid and HV grid of the area of Marseille-Nice (French Riviera).

%
\subsection{Convergence}
\label{subsec:Convergence}
%

Convergence is the main network analysis tool used in RTE. It is fully developped and maintained by RTE's R\&D teams.
Convergence is used for several time horizon usages: real time advance computations (state estimation, N-1 security analysis), operational planning and grid development. It embeds in a single platform (and single Graphical User Interface) many computation models: Load Flow computation model (named HADES), slow dynamics (ASTRE, to simulate voltage stability with online tap changers transformers actions), short circuit computation (COURCIRC), flow-based market coupling parameters, HVDC set point optimization (TOPAZE), several OPFs, and some others...
Convergence is also used in Coreso (European coordination center based in Brussels, Belgium) thanks to its coordination facilities: merging and analysis of European network files coming from D2CF (two days ahead capacity forecast), DACF (day ahead congestion forecast), IDCF (intra-day congestion forecast) processes. It is also used in real time in Coreso to merge and perform security analysis for European Snapshot files.
Concerning network data, Convergence is directly connected to the French national SCADA and the 7 regional SCADAs. Convergence gets every 5 minutes a full static network data from these 8 SCADAs, including equipment description, topology and state variables. Convergence performs merging of these 8 views of the French network (one national for Very High Voltage: 400kV and 225kV; 7 regional for High Voltage: 63kV to 150kV) to get a full consistent VHV+HV (63kV to 400kV) view of French transmission and sub-transmission Grid. Step-up transformers (20kV to 225kV or 400kV) are also included in grid data.
As already mentioned, in France Convergence is used with data coming directly from SCADAs. But it is also able to import and export grid data in UCTE and CIM\cite{CIMentsoe} formats. Eurostag format is also supported (import and export).
In this paper, 4 VHV+HV snapshot cases are published: they have more than 6000 nodes, more than 50\% are 63kV nodes.

%
\subsection{PEGASE}
\label{subsec:PEGASE}
%

The Pan European Grid Advanced Simulation and State Estimation (PEGASE) is a project of the 7\textsuperscript{th} Framework Program of the European Union~\cite{pegasesite}. Its goal was to develop new tools for the real-time control and operational planning of the pan-Euporean transmission network~\cite{pegase}. Specifically, new approaches were implemented for state estimation, dynamic security analysis, and steady state optimization. A dispatcher training simulator was also created.

The data accurately represents the size and complexity of the European high voltage transmission network. However, the data are fictitious and do not correspond to real world data. They can thus be used to validate methods and tools but should not be used for operation and planning of the European grid.

%
%
\section{Conversion of the data}
%
%

Some modifications to the original grid data were made in order to fit them into the \matpower{} format.

\smallskip

Concerning the snapshots, the generating units with pumped-storage capabilities could not be represented in \matpower{} format. Indeed, generating units with negative generation values are considered by \matpower{} as dispatchable loads. The way dispatchable loads and pumped-storage units are operated are significantly different, since pumped-storage units may work even if power prices are high because the stored energy is needed in the near future if power prices are expected to be higher. Moreover, reactive power production is not the same for dispatchable loads and for pumped-storage units. For these reasons, all generating units with Pmin lower than zero in our snapshot data have been converted into generating units with Pmin equals to zero.

Although some aspects of snapshots data might seem strange, that's the way they are in our network analysis tools. The most meshed 380kV bus is \textit{really} connected to 17 branches; this is not an artefact. Neighbouring countries are approximately represented with Ward equivalents; so these parts are not real. But they are really parts of data used in our SCADAs and network analysis tools. Some of these equivalent parts, as well as equivalent representation of 3-windings transformers, may lead to (small) negative reactances or resistances. To summarize, the data published in this article allow to build OPF problems that we really have to solve; iTesla data and PEGASE date share this property.

\smallskip

PEGASE data contains asymmetric shunt conductance and susceptance in the PI transmission line model of branches. However, \matpower{} format does not allow for asymmetry. As a result, we set the total line charging susceptance of branches to 0 per unit in the \matpower{} files. Instead, we used the nodal representation of shunt conductance and susceptance. This procedure leaves the power flow equations unchanged compared with the original PEGASE data. However, line flow constraints in the optimal power flow problem are modified.

%
%
\section{Description of the data}
%
%

In this section we give a few figures about the network cases.

Although pegase cases 89, 1354, 2869 and 9241 were already published in \matpower{} in 2015, we include them in description and result tables.

%
\subsection{General figures}
%

This first table gives for each case the number of buses, generating units, branches and transformers.

\smallskip

\noindent
\begin{tabular}{|l|r|r|r|r|}
\hline
Case Name&Bus.&Gen.&Bran.&Tran.\\
\hline
case89pegase&89&12&210&32\\
case1354pegase&1 354&260&1 991&234\\
case1888rte&1 888&297&2 531&405\\
case1951rte&1 951&391&2 596&486\\
case2848rte&2 848&547&3 776&558\\
case2868rte&2 868&599&3 808&606\\
case2869pegase&2 869&510&4 582&496\\
case6468rte&6 468&1 295&9 000&1 319\\
case6470rte&6 470&1 330&9 005&1 333\\
case6495rte&6 495&1 372&9 019&1 359\\
case6515rte&6 515&1 388&9 037&1 367\\
case9241pegase&9 241&1 445&16 049&1 319\\
case13659pegase&13 659&4 092&20 467&5 713\\
\hline
\end{tabular}

\smallskip

The next table gives the range of Voltage Levels that are included in each case.

\smallskip

\noindent
\begin{tabular}{|l|c|}
\hline
Case Name& Voltage Levels (kV)\\
\hline
case89pegase&380 220 150\\
case1354pegase&380 220\\
case1888rte&380 225 150 90 63 \& 24$\rightarrow$3\\
case1951rte&380 225 150 90 63 \& 24$\rightarrow$3\\
case2848rte&380 225 150 63 \& 45$\rightarrow$3\\
case2868rte&380 225 150 63 \& 45$\rightarrow$3\\
case2869pegase&380 220 150 110\\
case6468rte&380 225 150 90 63 \& 45$\rightarrow$3\\
case6470rte&380 225 150 90 63 \& 45$\rightarrow$3\\
case6495rte&380 225 150 90 63 \& 45$\rightarrow$3\\
case6515rte&380 225 150 90 63 \& 45$\rightarrow$3\\
case9241pegase&750 400 380 330 220 154 150 120 110\\
case13659pegase &750 400 380 330 220 154 150 120 ... \phantom{1} \\ & 110 \& 27$\rightarrow$0.4\\
\hline
\end{tabular}

\smallskip

The next table gives the number of buses for each main Voltage Level category.

\smallskip

\noindent
\hspace{-2ex}
\begin{tabular}{|l|r|r|r|r|r|}
\hline
&\multicolumn{5}{c|}{Number of Nodes per Voltage Level}\\
\hline
&\hspace{-1ex}$\ge$330kV& \hspace{-1ex}225kV & \hspace{-1ex}$\le$154kV & \hspace{-1ex}$\le$63kV & \hspace{-1ex}$\le$27k \\
Case Name
&               & \hspace{-1ex}220kV & $\hspace{-1ex}\ge$ 90kV  &  \hspace{-1ex}$\ge$45kV &         \\
\hline
case89pegase  & 50 & 5 & 34 & 0 & 0 \\
case1354pegase  & 241 & 1113 & 0 & 0 & 0 \\
case1888rte  & 349 & 1174 & 61 & 8 & 296 \\
case1951rte  & 350 & 1185 & 62 & 8 & 346 \\
case2848rte  & 347 & 1177 & 59 & 915 & 350 \\
case2868rte  & 351 & 1193 & 59 & 918 & 347 \\
case2869pegase  & 629 & 1748 & 492 & 0 & 0 \\
case6468rte  & 524 & 1274 & 1183 & 3151 & 336 \\
case6470rte  & 525 & 1277 & 1183 & 3150 & 335 \\
case6495rte  & 525 & 1277 & 1184 & 3152 & 357 \\
case6515rte  & 525 & 1283 & 1184 & 3153 & 370 \\
case9241pegase   & 1945 & 3185 & 4111 & 0 & 0 \\
\hspace{-1ex}case13659pegase  & 1945 & 3185 & 4111 & 0 & 4418 \\
\hline
\end{tabular}

\smallskip

The next table gives the number of branches with negative resistance R and the number of branches with negative reactance X. There is no branch with both negative R and X.

\smallskip

\noindent
\begin{tabular}{|l|r|r|}
\hline
Case Name & Branches R$<$0 & Branches X$<$0 \\
\hline
case89pegase & 0     & 0 \\
case1354pegase & 0     & 0 \\
case1888rte & 0     & 77 \\
case1951rte & 0     & 76 \\
case2848rte & 0     & 75 \\
case2868rte & 0     & 77 \\
case2869pegase & 0     & 0 \\
case6468rte & 0     & 80 \\
case6470rte & 0     & 80 \\
case6495rte & 0     & 80 \\
case6515rte & 0     & 80 \\
case9241pegase & 75    & 16 \\
case13659pegase & 78    & 16 \\
\hline
\end{tabular}%

%
\subsection{Impedances and voltages}
%

In this section we give an outlook on data, using illustrative graphs.
For 3 cases, a first graph shows norm of impedances of all lines, in descending order and logarithmic scale.
A second graph shows the complex values of voltages of all buses; the circles on these graphs are corresponding to 0.9 and 1.1 p.u.

Among the 4 iTesla cases, we show these two graphs for case2868rte.
Among the 4 RTE snapshot cases, we show these two graphs for case6515rte.
For PEGASE cases, we chose to show the graphs for the largest case: case13659pegase.

\begin{figure}[H]
  \centering
    \includegraphics[width=.5\textwidth]{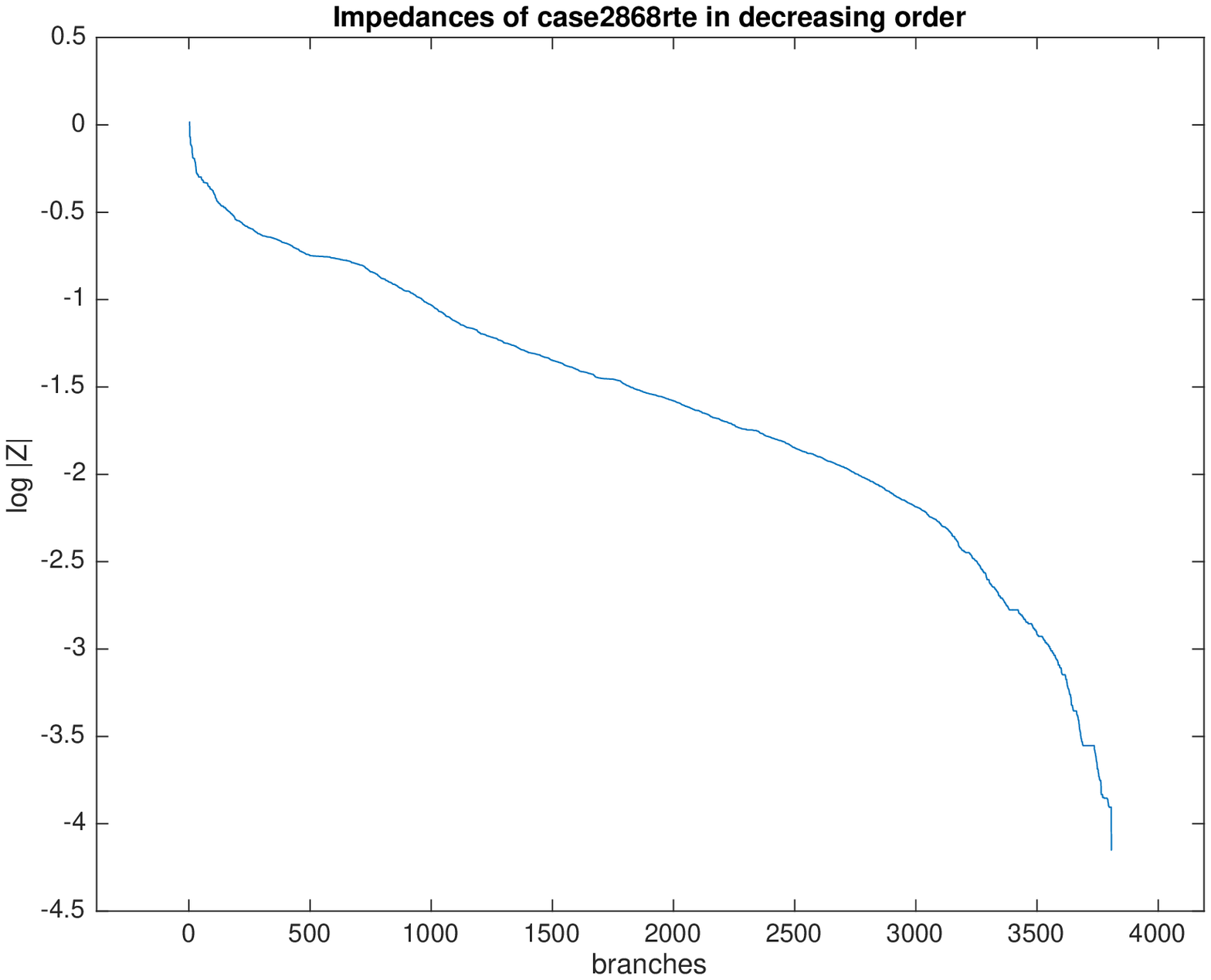}
  \label{fig:case2868impedance}
\end{figure}
\begin{figure}[H]
  \centering
    \includegraphics[width=.5\textwidth]{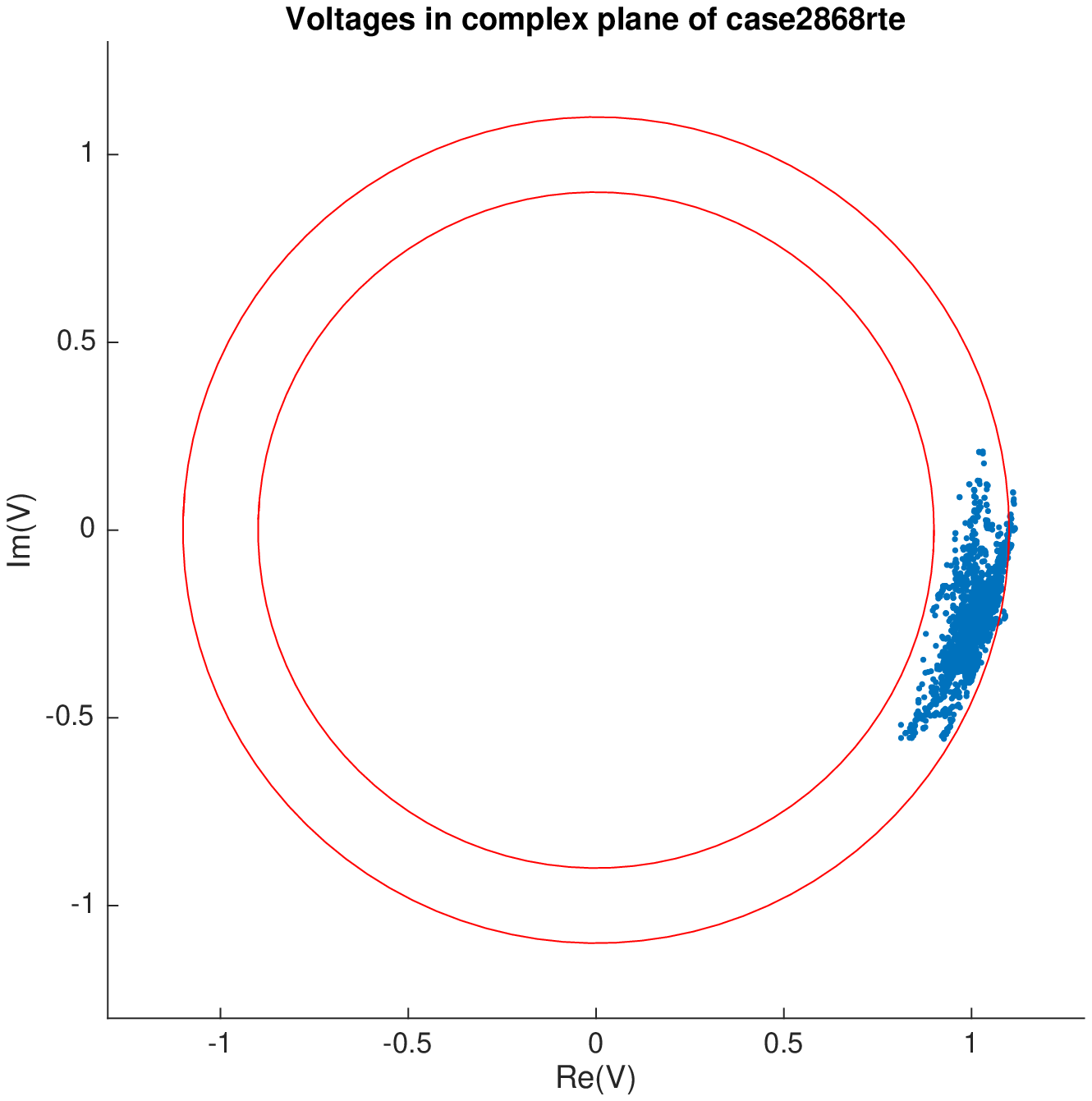}
  \label{fig:case2868voltage}
\end{figure}

\begin{figure}[H]
  \centering
    \includegraphics[width=.5\textwidth]{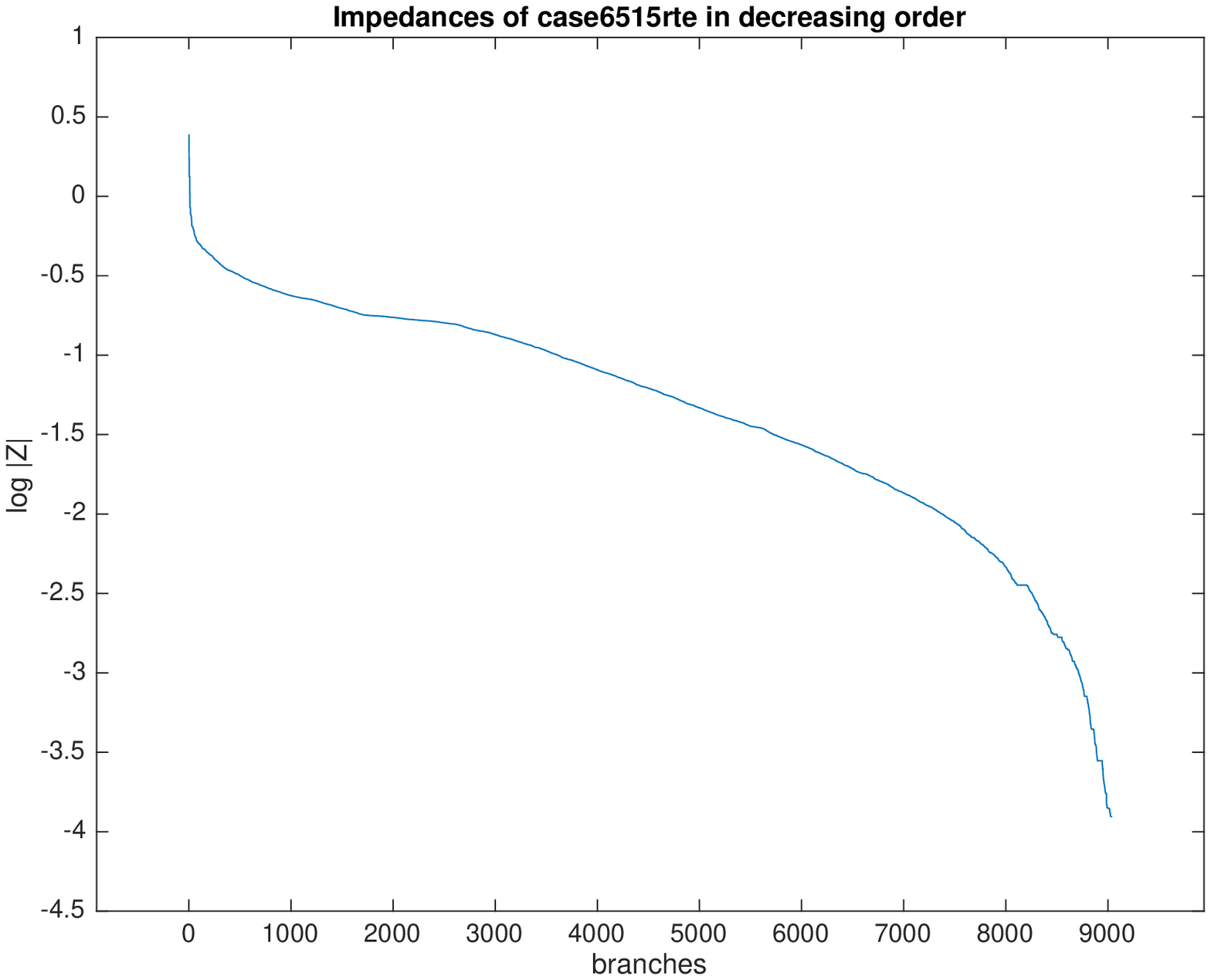}
  \label{fig:case6515impedance}
\end{figure}
\begin{figure}[H]
  \centering
    \includegraphics[width=.5\textwidth]{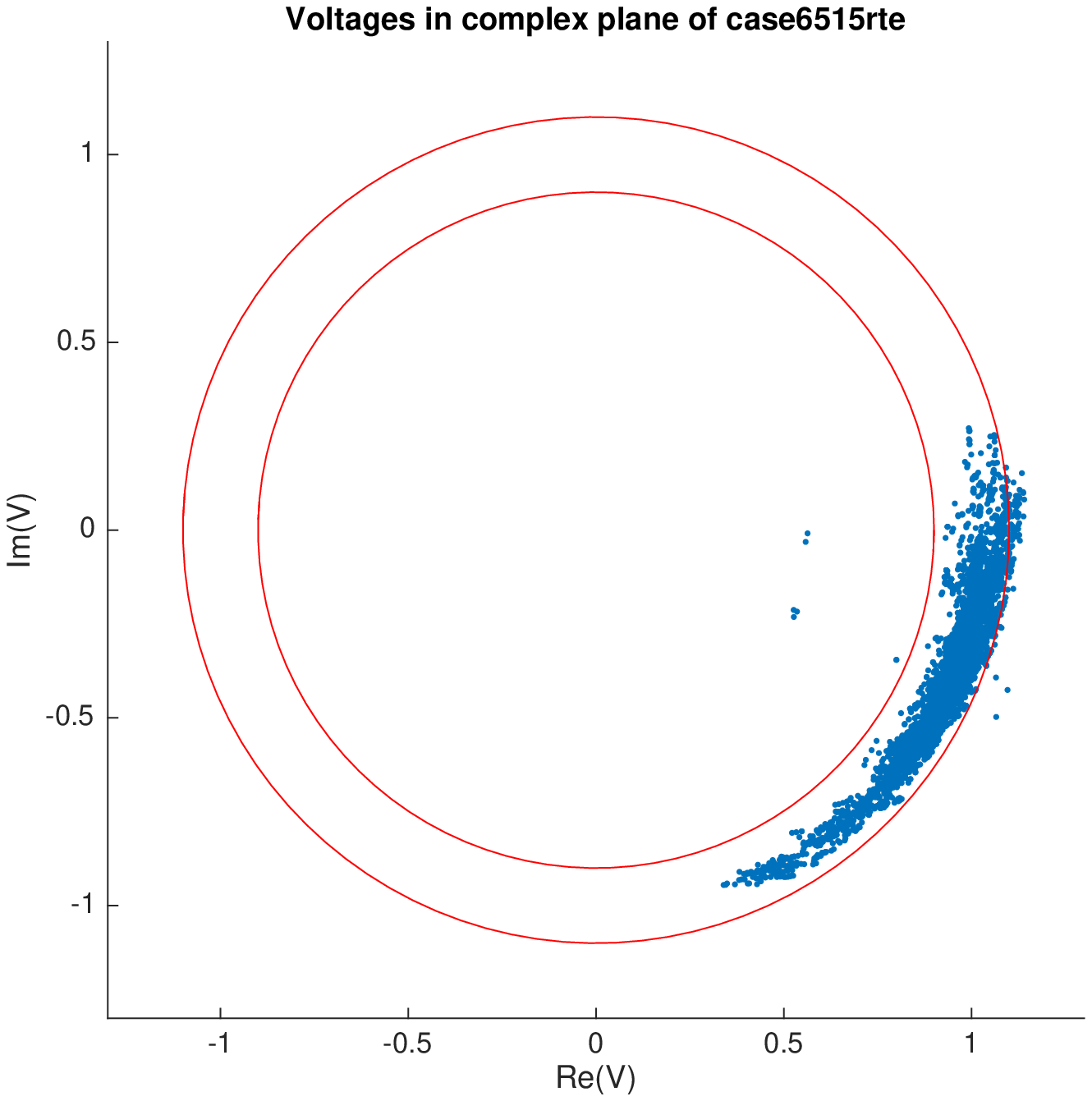}
  \label{fig:case6515voltage}
\end{figure}

\begin{figure}[H]
  \centering
    \includegraphics[width=.5\textwidth]{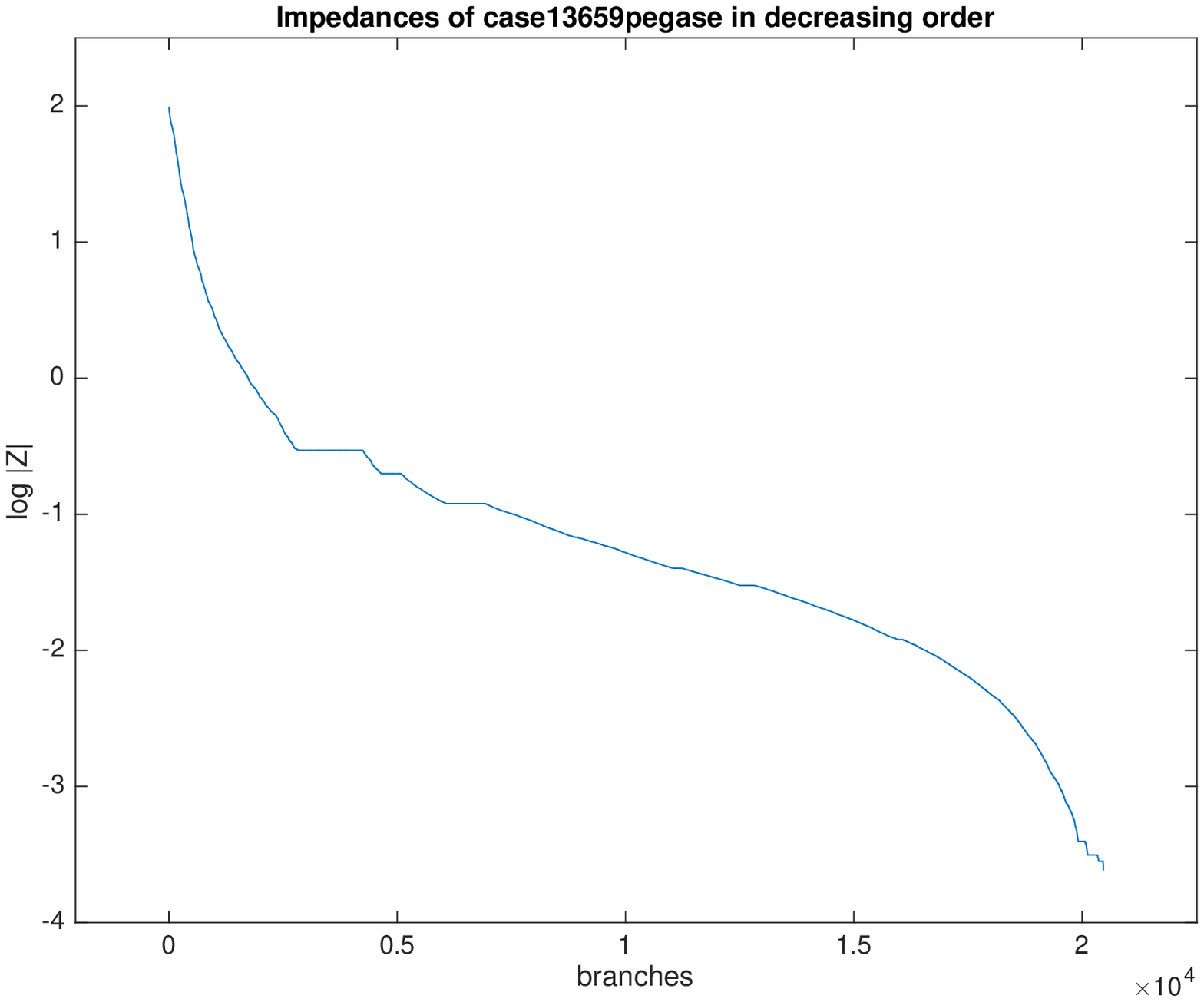}
  \label{fig:case13659impedance}
\end{figure}

\begin{figure}[H]
  \centering
    \includegraphics[width=.5\textwidth]{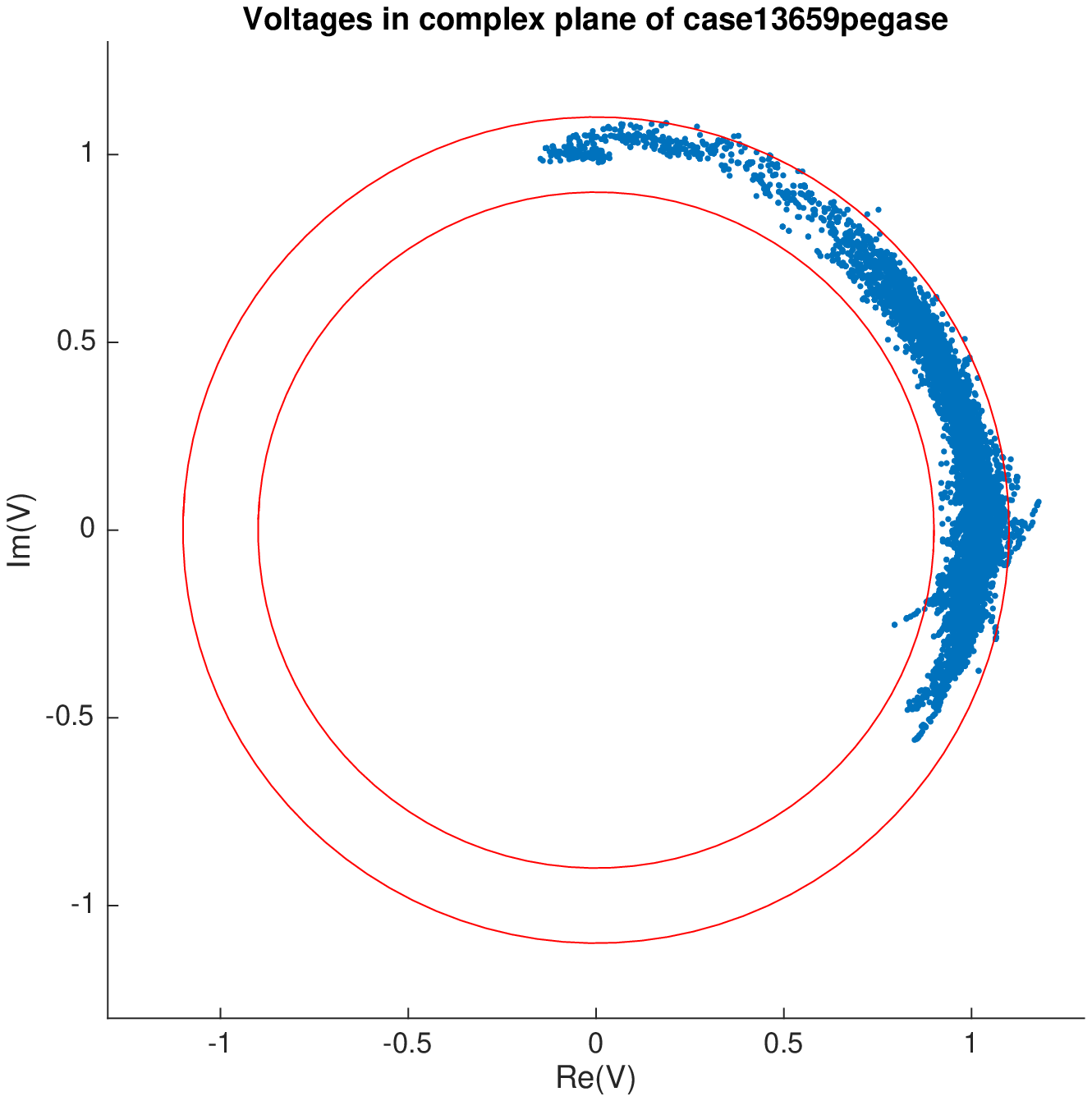}
  \label{fig:case13659voltage}
\end{figure}

\section{Numerical results}
%
%

%
\subsection{Using Knitro to find feasible solutions}
%

With Knitro~\cite{knitro} as interior point solver, we were able to run \matpower{} Optimal Power Flow ('runopf').
We found out that \matpower{} option for Knitro 'xtol' is set to 0.0001, which is too large. Using an option file for Knitro (file 'knitro.opt'),  we could modify 'xtol' parameter of Knitro.

\smallskip

\begin{tabular}{|l|c|}
\hline
\multicolumn{2}{|c|}{Options for Knitro in knitro.opt} \\
\hline
xtol   & 1e-8 \\
feastol & 5e-6 (default value of \matpower{}) \\
opttol & 1e-4 (default value of \matpower{}) \\
maxit &  1000 \\
\hline
\end{tabular}

\smallskip

All OPF could be solved with this parameterization of Knitro.

Using Knitro to solve a non convex problem gives a solution which is expected to satisfy local optimality conditions. In this paper, the fact that the result of Knitro is a local solution is not used. Instead, we simply use the property that Knitro's solution is a feasible solution for OPF problem. Any other software, any other optimization method, any other heuristic might be used to find better feasible solution of OPF.

For instance, a method to find a better feasible solution of OPF could be:
\begin{enumerate}
\item From original OPF, create a modified OPF with better numerical properties (e.g. agregation of electrical nodes with very small electric distances as in \cite{cdc2015bis}); eventually also add penalty terms in OPF \cite{lavaei_mesh,lavaei_allerton2014,cdc2015bis,molzahn_josz_hiskens_panciatici-Laplacian_Objective}.
\item Solve to global optimality the modified OPF.
\item Use the solution of the modified OPF to build a feasible solution for original OPF.
\end{enumerate}
In any case, finding better feasible solutions for OPF will not be enough: finding better lower bounds will also help to achieve global optimality for OPF.

%
\subsection{Basic evaluation of optimality gap}
%

Semidefinite programming is now known to be able to give good lower bounds for OPF problem. However, comparisons with basic lower bounds may not be avoided.

In the next table, the first column contains result of DCOPF with Mips (default solver of \matpower{}), without flow limits. As our cost functions are all linear with costs equal to 1, it means that our objective function is simply minimizing the total sum of generation, or equivalently losses minimization. As losses are neglected in DC modelling, optimal value of our DCOPFs is to be equal to the sum of loads. So why writing here the results of these DCOPF ? Because it is a trivial lower bound of ACOPF. In case cost functions really embed generating units different costs, results of DCOPF without flow constraints will also be a basic lower bound of ACOPF. In the next table, the second column shows result of ACOPF (without flow constraints) solved with Knitro. Using the basic lower bound computed with DCOPF, optimality may be computed in the third column. 

Case9241pegase and case13659pegase have branches with negative resistances; for these two cases, DCOPF result is not a strict lower bound of OPF, since pushing a lot of power through branches with negative resistances would create active power. For this reason, the two last optimality values are not sure.

\smallskip

\noindent \hspace{-1ex}
\begin{tabular}{|l|r|r|r|}
\hline
\multicolumn{4}{|c|}{Losses minimization without flow constraints} \\
\hline
\multicolumn{1}{|r|}{Algorithm:}  & \multicolumn{1}{c|}{DCOPF} & \multicolumn{1}{c|}{OPF}   & \multirow{2}[0]{*}{Optimality} \\
\multicolumn{1}{|r|}{Solver:}       & \multicolumn{1}{c|}{Mips}  & \multicolumn{1}{c|}{Knitro} &  \\
\hline
    case89pegase & 5 733.4 & 5 817.6 & 1.47\% \\
    case1354pegase & 73 059.7 & 74 060.4 & 1.37\% \\
    case1888rte & 59 110.5 & 59 769.9 & 1.12\% \\
    case1951rte & 80 656.5 & 81 724.2 & 1.32\% \\
    case2848rte & 52 562.3 & 53 020.9 & 0.87\% \\
    case2868rte & 78 826.3 & 79 783.4 & 1.21\% \\
    case2869pegase & 132 447.2 & 133 980.7 & 1.16\% \\
    case6468rte & 85 296.9 & 86 791.8 & 1.75\% \\
    case6470rte & 96 592.4 & 98 308.0 & 1.78\% \\
    case6495rte & 103 916.1 & 105 943.6 & 1.95\% \\
    case6515rte & 107 264.0 & 109 561.2 & 2.14\% \\
    case9241pegase & 312 411.0 & 315 888.5 & 1.11\%? \\
    case13659pegase & 381 773.4 & 386 107.5 & 1.14\%? \\
\hline
\end{tabular}

\smallskip

Although solving ACOPF with Knitro usually gives only a local optimum, without any information whether this optimum might be global, interestingly we can see that all Knitro (local) solutions are 0.87\% to 2.14\% optimal. To our opinion, any method that gives results with larger proven optimality bounds can not claim to give good results. Moreover, precision of numerical methods should always be compared to the precision of the pair (DCOPF;easily obtained feasible solution).

%
\subsection{OPF with or without flow constraints}
%

We made three series of computation. First series was with flow limits in terms of apparent power (in MVA, using option 'S' of \matpower{}), second series was with flow limits in terms of current (in Amperas, using option 'I' of \matpower{}). Third series was without any limit on flows. Results are shown in the next table.

\smallskip

\noindent \hspace{-2ex}
\begin{tabular}{|l|r|r|r|}
\hline
\multicolumn{4}{|c|}{OPF for losses minimization with/without flow limits} \\
\hline
Case Name & \multicolumn{1}{c|}{\hspace{-1ex}Flow Lim. 'S'} & \multicolumn{1}{c|}{\hspace{-1ex}Flow Lim. 'I'} & \multicolumn{1}{c|}{\hspace{-1ex}No Limit} \\
\hline
case89pegase & 5 819.8 & 5 817.6 & 5 817.6 \\
case1354pegase & 74 069.4 & 74 064.2 & 74 060.4 \\
case1888rte & 59 805.1 & 59 808.5 & 59 769.9 \\
case1951rte & 81 737.7 & 81 737.4 & 81 724.2 \\
case2848rte & 53 021.8 & 53 021.9 & 53 020.9 \\
case2868rte & 79 794.7 & 79 794.5 & 79 783.4 \\
case2869pegase & 133 999.3 & 133 993.5 & 133 980.7 \\
case6468rte & 86 860.0 & 86 841.8 & 86 791.8 \\
case6470rte & 98 345.5 & 98 325.4 & 98 308.0 \\
case6495rte & 106 283.4 & 106 215.7 & 105 943.6 \\
case6515rte & 109 804.2 & 109 767.8 & 109 561.2 \\
case9241pegase & 315 912.7 & 315 903.3 & 315 888.5 \\
\hspace{-1ex}case13659pegase& 386 107.5 & 386 107.5 & 386 107.5 \\
\hline
\end{tabular}

\smallskip

 Reader will immediately notice that values with/without flow limits are very similar, and even identical for the largest case (case13659pegase). This is not abnormal. RTE snapshots are real observed data, Pegase and iTesla data were constructed to be realistic. In real life grid operation, N-1 rule implies that almost everywhere, lines are within their limits. Moreover, French grid was historically built to optimally serve load and production schemes for all situations; policy for construction of new uncontrollable renewable energy production (wind, solar) was to encourage small units distributed throughout the whole country, with limited impact on transmission grid so far. This situation is evolving fast, but snapshot data are from year 2013.
Occasionnaly, some lines may be over their limits: temporary admissible limits are used by operators, they usually have a 20 minutes time window to take corrective actions.

%
%

%
\subsection{Lower bounds with SDPOPF}
%

In this section we  tried to obtain lower bounds via semidefinite programming. We used SDPOPF solver provided by Daniel K. Molzahn \cite{dan2013} in \matpower{}, with Sedumi \cite{sturm-1999,sedumigithub} and Mosek \cite{mosek71} SDP solvers.

We would like to address a special thank to Dan K. Molzahn for helping us using his SDPOPF solver.
He pointed out to us that in order to have lower SDP bounds of our original OPF problems, some options of SDPOPF has to be set:

\smallskip

\noindent \hspace{-2ex}
\begin{tabular}{|l|}
\hline
\multicolumn{1}{|c|}{Options for SDPOPF in order to compute lower bounds} \\
\hline
 \texttt{mpopt.sdp\_pf.eps\_r = -inf} \\ 
$\rightarrow$ Do not enforce a minimum resistance \\
 \texttt{mpopt.sdp\_pf.min\_Pgen\_diff = 0} \\ 
$\rightarrow$ Disable enforcing fixed value when small range for Pgen \\
 \texttt{mpopt.sdp\_pf.min\_Qgen\_diff = 0} \\ 
$\rightarrow$ Disable enforcing fixed value when small range for Qgen \\
 \texttt{mpopt.sdp\_pf.max\_line\_limit = inf} \\ 
$\rightarrow$ Disable elimination of large line limits \\
 \texttt{mpopt.sdp\_pf.max\_gen\_limit = inf} \\ 
$\rightarrow$ Disable elimination of large generation limits \\
\hline
\end{tabular}

\smallskip

In the next table we compare the computed bounds; we expect to have each line sorted in ascending order: 
\begin{enumerate}
\item DCOPF value, basic lower bound of the optimal value, 
\item then the value of the SDP relaxation computed by SDPOPF, ideally the values obtained by the two different SDP solvers would be equal, 
\item and then, in the OPF column, a feasible solution of our OPF problem, giving an upper bound to the OPF problem.
\end{enumerate}

\smallskip

\noindent \hspace{-2ex}
\begin{tabular}{|l|r|r|r|r|}
\hline
\multicolumn{5}{|c|}{Losses minimization without flow constraints} \\
\hline
\multicolumn{1}{|r|}{Algorithm:}  & \multicolumn{1}{c|}{DCOPF} & \multicolumn{1}{c|}{SDPOPF} & \multicolumn{1}{c|}{SDPOPF} & \multicolumn{1}{c|}{OPF} \\
\multicolumn{1}{|r|}{Solver:}  & \multicolumn{1}{c|}{Mips} & \multicolumn{1}{c|}{Sedumi} & \multicolumn{1}{c|}{Mosek} & \multicolumn{1}{c|}{Knitro} \\
\hline
case89pegase & 5 733.4 & 5 817.6 & 5 817.6 & 5 817.6 \\
case1354pegase & 73 059.7 & 74 052.8 & 74 049.5 & 74 060.4 \\
case1888rte & 59 110.5 & 59 572.0 & 59 557.7 & 59 769.9 \\
case1951rte & 80 656.5 & 81 718.7 & 81 706.4 & 81 724.2 \\
case2848rte & 52 562.3 & 53 006.6 & 52 986.4 & 53 020.9 \\
case2868rte & 78 826.3 & 79 782.9 & 79 769.1 & 79 783.4 \\
case2869pegase & 132 447.2 & 133 970.9 & 133 964.6 & 133 980.7 \\
case6468rte & 85 296.9 & 86 754.5 & 86 726.2 & 86 791.8 \\
case6470rte & 96 592.4 & 98 305.0 & 98 277.0 & 98 308.0 \\
case6495rte & 103 916.1 & 105 969.7 & 105 919.4 & 105 943.6 \\
case6515rte & 107 264.0 & 109 560.7 & 109 533.5 & 109 561.2 \\
case9241pegase & 312 411.0 & 310 723.5 & 310 697.1 & 315 888.5 \\
\hspace{-1ex}case13659pegase & 381 773.4 & 381 047.8 & 381 027.7 & 386 107.5 \\
\hline
\end{tabular}

\medskip

First, it is very interesting to note that extra large problems (6 to 13 thousends of buses) could be addressed by SDPOPF without reaching the limits of our 48GB RAM computer. We also tried SDPT3 instead of Sedumi or Mosek, but we encountered matlab exceptions with cases larger than 6000 buses.

Second, looking more precisely into log files, we could see that in \textit{all} cases, Sedumi ended its computation with message
\texttt{Run into numerical problems} and Mosek with message \texttt{Mosek error: MSK\_RES\_TRM\_STALL()}. It means that SDP solvers are not totally mature yet to solve our problems, and we have to take care using their results.
For instance, case6494rte has its Sedumi SDP lower bound larger than OPF/Knitro value; obviously the SDP value is not well computed. For this reason, we think that results obtained with Sedumi are not precise enough to claim that, e.g., case6515rte would be solved to global optimality with an error smaller than 0.5MW.

Third, SDP lower bounds are smaller than DCOPF values for the two largest cases. This is certainly due to the presence of negative resistances in these two cases: DCOPF values are not lower bounds of OPF and the only lower bounds we have are the SDP ones. This point was already mentioned in \cite{josz-molzahn-2015}, when analyzing results of the Shor relaxation in Table 2 of \cite{josz-molzahn-2015}.

Last, with Mosek all results seems consistent. For the two largest cases, global optimality is proven only to 1.7\% and 1.3\%, which are not very good values (order of magnitude of the losses). case1888rte' optimality is 0.36

Note that OPF (Knitro) results are not impacted by precision issues: running again the OPF (Knitro) column with feastol equals to 1e-10 gives the same results.

\smallskip

As a conclusion, in the next table, we compare global optimality proofs obtained with DCOPF (Mips solver) and SDPOPF (Mosek solver).

\smallskip

\noindent
\begin{tabular}{|l|r|r|}
\hline
\multicolumn{3}{|c|}{Global Optimality proofs} \\
\hline
      & \multicolumn{1}{c|}{DCOPF} & \multicolumn{1}{c|}{SDPOPF} \\
\hline
case89pegase & 1.47\% & 0.00\% \\
case1354pegase & 1.37\% & 0.01\% \\
case1888rte & 1.12\% & 0.36\% \\
case1951rte & 1.32\% & 0.02\% \\
case2848rte & 0.87\% & 0.07\% \\
case2868rte & 1.21\% & 0.02\% \\
case2869pegase & 1.16\% & 0.01\% \\
case6468rte & 1.75\% & 0.08\% \\
case6470rte & 1.78\% & 0.03\% \\
case6495rte & 1.95\% & 0.02\% \\
case6515rte & 2.14\% & 0.03\% \\
case9241pegase & not valid & 1.67\% \\
case13659pegase & not valid & 1.33\% \\
\hline
\end{tabular}

%
\subsection{Global optimality quest}
%

Our goal is to be able, in the near future, to prove global optimality (with precision $10^{-5}$ or $10^{-6}$) for all these cases.

From operational point of view, it is worth to spend time to compute a \textit{global} optimum? Non global method such as interior point methods have been sucessfully used for about 20 years to solve OPF problems, so is \textit{global} optimality necessary?
Isn't it only a game for scientists, with no industrial consequence?

Our answer to this question is in a larger view of optimization methods for grid operations and development.
Once continuous OPF will be solved to global optimality (for losses minimization, but also with all kinds of generation costs), we'll start addressing global optimization of OPF with discrete variables (e.g. on/off generating units statuses, but also PST taps, discrete shunts, topology choices...).
When global optimality of OPF will be easily available, we'll be able to solve bilevel programs such as in \cite{Mitsos2009}, for which global optimality of OPF is necessary as OPF are subproblems of a wider framework.

%
%
\section{Mathematical format}
%
%

In addition to the grid data, we provide a code \texttt{qcqp\_opf.m} that converts  any MATPOWER test case data into a standard mathematical optimization format. Its purpose is to allow members of the applied mathematics community to evaluated their methods on the test cases without requiring any knowledge in power systems.  Precisely, the MATPOWER test cases are converted into large-scale sparse quadratically-constrained quadratic programs (QCQP). Indeed, the optimal power flow problem can be viewed as an instance of quadratically-constrained quadratic programming. In order for this to be true, we consider the objective function of the optimal power flow problem to be a linear function of active power. Higher degree terms are discarded from the objective function. Moreover, current line flow constraints are enforced instead of apparent line flow constraints in order to have quadratic constraints only. The optimal power flow problem remains non-convex and non-deterministic polynomial-time hard despite the slightly simplified framework we consider. Notice that for rte and pegase cases (minimization of the total generation) there is no simplification: QCQP formulation is equivalent to OPF with current line flow constraints. 

The standard format we use is described below where $x$ is a column vector of size nVAR:\\\\
\textbf{QCQP:}
\begin{equation}
\inf_x ~~~ x' C x + c ~~~
\end{equation}
subject to nEQ equality constraints
\begin{equation}
x' A_k x = a_k~,~~~ \forall k =1... \text{nEQ},
\end{equation}
and subject to nINEQ inequality constraints
\begin{equation}
x' B_k x \leqslant b_k~,~~~ \forall k =1...\text{nINEQ},
\end{equation}
where $C$, $A_k$'s, and $B_k$'s are squares matrices of size nVAR, $a,b$ are column vectors, and the apostrophe stands for conjugate transpose.
The code provides matrices that either complex, Hermitian, or real symmetric depending on an input parameter (see comments in code for details). Depending on this input parameter, column vectors $x,a$ and $b$ are either complex or real.

The following table shows the size of the QCQP instances in real numbers. The right column corresponds to the percentage of monomials that have a non-zero coefficient in the objective or constraints compared to the total number of possible monomials (in the case of a fully dense QCQP problem). It shows that the OPF is a very sparse problem so we believe that it is possible to solve these instances to global optimality.

\medskip

\noindent
\begin{tabular}{|l|r|r|r|r|r|}
\hline
Case Name&nVAR&nEQ&nINEQ & Spa. (\%)\\
\hline
case89pegase&178&154 & 380 & 5.23\\
case1354pegase&2 708 & 2 188 & 6 612 & 0.19 \\
case1888rte&3 776 & 3 222 & 9 036 & 0.13 \\
case1951rte&3 902 & 3 162 & 9 580 & 0.12 \\
case2848rte&5 696 & 4 904 & 11 742 & 0.08 \\
case2868rte&5 736 & 4 850 & 12 070 & 0.08 \\
case2869pegase&5 738 & 4 718 & 13 264 & 0.10 \\
case6468rte&12 936 & 11 652 & 20 130 & 0.04 \\
case6470rte&12 940 & 11 588 & 21 864 & 0.04 \\
case6495rte&12 990 & 11 562 & 22 064 & 0.04 \\
case6515rte&13 030 & 11 576 & 22 200 & 0.04 \\
case9241pegase&18 482 & 15 592 & 36 852 & 0.03 \\
case13659pegase& 27 318 & 19 134 & 43 686 & 0.02 \\
\hline
\end{tabular}

%
%
\section{Conclusion}
%
%

In this paper our goal is to publish very realistic data, being used every day by a large Transmission System Operator. Elementary description of data and their origin are included. Preliminary OPF results are also provided.

We aim to publish new versions of this document in the future, with additional numerical results (e.g. better upper or lower bounds) coming either from our own research activities or from other public academic works.

All cases mentioned in this paper are included in this arXiv publication (in the source tar file that can be donwloaded from arXiv), except the PEGASE cases that were already published in \matpower{} in 2015.

A MATLAB code to transform OPF data to standard QCQP mathematical optimization format is also included, with the hope that it will help mathematicians address these problems without power system skills.

Moreover, a few m-files are also included in the tar file, in order to help scientific OPF community to uses these data and to reproduce results.

%
%
%
%
\bibliography{mybib}{}
\bibliographystyle{plain}

\end{document}